\documentclass[12pt,a4paper]{article}
\usepackage{template}

\begin{document}
\pagestyle{plain}
\title{A Galerkin least-square stabilisation technique for hyperelastic biphasic soft tissue}
\date{}
\author{Julien Vignollet}
\author{Chris J. Pearce}
\author{Lukasz Kaczmarczyk}
\affil{School of Engineering, Glasgow University}
\maketitle

\begin{abstract}
An hyperelastic biphasic model is presented. For slow-draining problems (permeability less than 1$\times$10$^{-2}$ mm$^{4}$N$^{-1}$s$^{-1}$), numerical instabilities in the form of non-physical oscillations in the pressure field are observed in 3D problems using tetrahedral Taylor-Hood finite elements. As an alternative to considerable mesh refinement, a Galerkin least-square stabilization framework is proposed. This technique drastically reduces the pressure discrepancies and prevents these oscillations from propagating towards the centre of the medium. The performance and robustness of this technique are demonstrated on a 3D numerical example.
\end{abstract}


\section{Introduction}
The theory of porous media (TPM) is a powerful and yet simple tool to model multi-phase media, primarily comprising a solid and a fluid (which are referred to as "constituents"). Due to its modularity, this theory also represents a practical framework to include other effects such as mass exchange, chemical reactions and electrochemical phenomena. There is no restriction as to which material models can be used, e.g. anisotropy or viscosity and is therefore well suited for the analysis of hydrogels, polymeric foams and hydrated biological soft tissues (e.g. cartilage and intervertebral disc).

\begin{figure}[ht]
\begin{centering}
\includegraphics[width=0.8\textwidth]{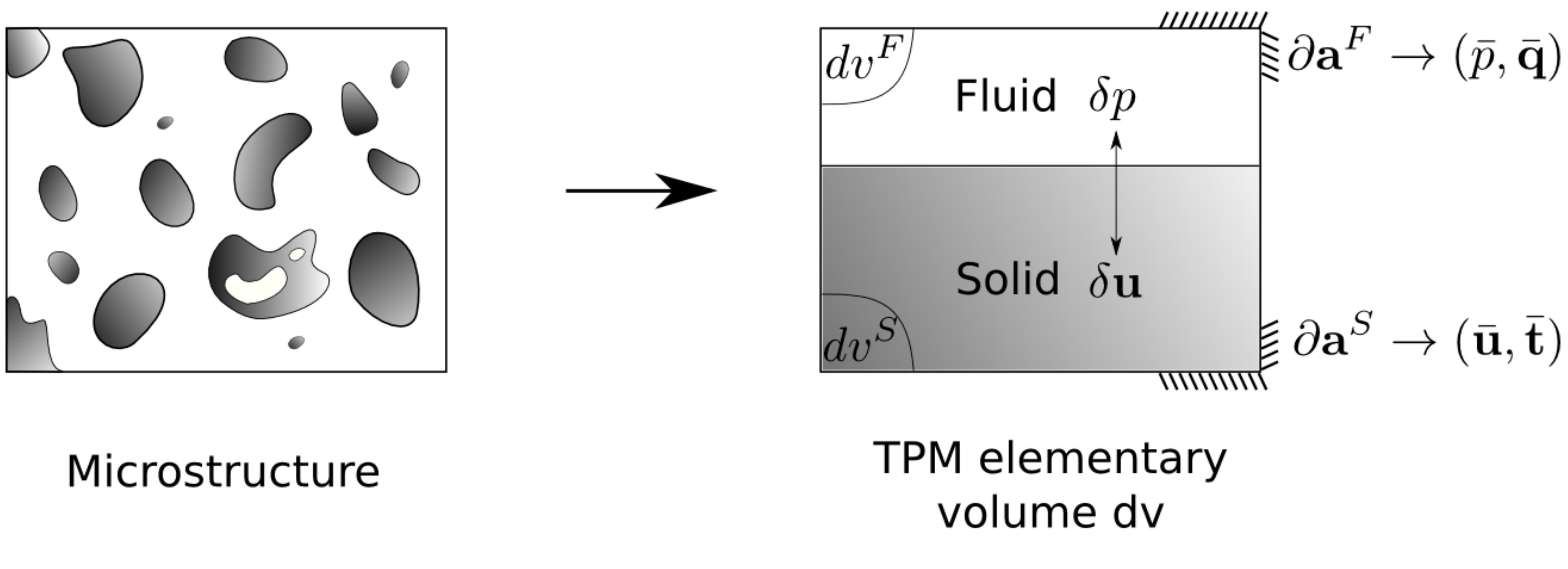}
\par\end{centering}
\caption{The TPM representation (right) of the microstructure (left); at each material point, fluid and solid coexist. Displacement and pressure fields are coupled. \label{fig:mixture}}
\end{figure}

The concepts behind the TPM find their roots in diffusion and soil mechanics problems formulated in the nineteenth century (for historical developments, see \cite{Atkin1976} and \cite{deBoer2005}). The TPM is a homogenized macroscopic representation of the porous media. It is a continuum-based model which fully couples the fluid and the solid (see Fig. \ref{fig:mixture}). It overcomes the difficulty of obtaining an accurate geometrical description of the microstructure by using the concept of volume fractions to ``smear'' the constituent properties over a control space to obtain properties of the overall mixture. This is described in the following section.

The authors have employed TPM in the development of an hyperelastic biphasic swelling model for modelling the intervertebral disc. This application illustrates slow-draining problems, where a porous medium with low permeability (in the range 1$\times$10$^{-2}$$\--$1$\times$10$^{-3}$ mm$^{4}$N$^{-1}$s$^{-1}$) is rapidly loaded. Such low permeabilities are typical for soft tissues (e.g. cartilage, brain tissue), mortar or homogeneous clays (\cite{Ateshian1997}, \cite{Loosveldt2002}, \cite{Terzahi1996}).

This model has been implemented in a finite element framework, employing Taylor-Hood (quadratic shape functions for the solid displacement, linear shape functions for the pressure) tetrahedral elements, aiming to fulfil the \emph{inf-sup condition} (see \cite{Chapelle1993} and \cite{Brezzi1991}). However, numerical instabilities manifest in the form of non-physical oscillations in the pressure field, which is a shortcoming already observed in the past (see for example \cite{Stokes2010} for a recent review of biological applications). Vermeer and Verruijt \cite{Vermeer1981} explain that these instabilities occur because loads applied to free-flow boundary conditions may lead to singularities in the derivatives of the pressure field. They also derive a lower bound critical time-step for one-dimensional problems, suggesting the requirement for large time-steps to overcome this issue, often incompatible with fast loading rates.

Several stabilisation techniques have been proposed in the context of Biot's consolidation problems for small deformations (e.g. \cite{Korsawe2005} using least-squares mixed finite element methods, and \cite{Aguilar2008} by perturbation of the flow equation). The current work proposes to stabilize the pressure oscillations in the context of TPM for finite deformation problems, using a Galerkin Least-Square (GLS) formulation based on \cite{Truty2001}. In order to focus only on the stabilization aspect, only a biphasic mixture is considered, that is a porous medium composed of two constituents $\alpha$: an isotropic, non viscous hyperelastic solid ($\alpha=S$) and an ideal fluid ($\alpha=F$).

\section{The biphasic model}
The biphasic model presented in this section is based on \cite{Boer1998} and \cite{Ehlers2002a}. A few preliminary assumptions are made to keep the derivation as simple as possible in order to focus on the stabilization aspect. First, the quasi-static problem of small biological tissues is herein considered, thus neglecting external body forces. Second, the constituents are assumed immiscible and no density supply is allowed. Third, it is assumed that the whole space is occupied by either of the constituents. Finally, intrinsic incompressibility is assumed for both constituents. It is important to note that the last assumption is not equivalent to incompressibility of the whole mixture.

The most fundamental concept of mixture theory (established as early as in \cite{Stefan1871}) asserts that at any time $t$ and at each spatial point $\mathbf{x}$ of the continuum, particles of both constituents $\alpha$ coexist. This implies that any elementary volume $dv$ is simultaneously occupied by both phases and is split into partial elementary volumes $dv^\alpha$. The volume fractions can then be defined as:

\begin{equation}
n^{\alpha} \left( \mathbf{x} , t \right)= \frac {dv^{\alpha}} {dv} \qquad \alpha = \left\{ S,F \right\}
\label{eq:volume fraction}
\end{equation}

Assuming that there is no gas in the mixture, the saturation condition is:

\begin{equation}
	n^{S} + n^{F} = 1 
\label{eq:saturation condition}
\end{equation}

The useful relationship between the apparent density $\rho^{\alpha}$ and true density $\rho_{T}^{\alpha}$ of each constituent is defined as follows:

\begin{equation}
\rho^{\alpha}=\frac{dm^{\alpha}}{dv}=\frac{dm^{\alpha}}{dv^{\alpha}}\frac{dv^{\alpha}}{dv}=\rho_{T}^{\alpha}n^{\alpha}
\label{eq:connect density to true density}
\end{equation}

The other main principles of mixture theory are summarised in \emph{``Truesdell's metaphysical principles''} \cite{Truesdell1984}
	\footnote{``Truesdell's metaphysical principles'': 1) All properties of the mixture must be mathematical consequences of properties of the constituents; 2) So as to describe the motion of a constituent, we may in imagination isolate it from the rest of the mixture, provided we allow properly for the actions of the other constituents upon it; 3) The motion of the mixture is governed by the same equations as is a single body.}. 
The first consequence of these principles is that each constituent ${\alpha}$ of the mixture is described by its own state of motion $\bchi^{\alpha}$ (see Eq. \ref{eq:motion}), relating at time $t$ the position vector $\mathbf{X}^{\alpha}$ of a particle in the reference configuration to its position in the current configuration $\mathbf{x}^{\alpha}$.

\begin{equation}
\mathbf{x}^{\alpha} =  \bchi^{\alpha} \left( \mathbf{X^{\alpha}} , t \right)
\label{eq:motion}
\end{equation}

This subsequently implies the existence of a velocity field $\mathbf{v}^{\alpha}$ and a deformation gradient $\mathbf{F}^{\alpha}$ for each constituent.

It proves convenient to describe the fluid velocity as a velocity relative to the solid constituent, leading to the introduction of the seepage velocity $\mathbf{w}$:

\begin{equation}
\mathbf{w}=n^{F}\left(\mathbf{v}^{S}-\mathbf{v}^{F}\right)
\label{eq:seepage velocity}
\end{equation}

As mentioned earlier, each constituent's intrinsic incompressibility does not imply mixture incompressibility. Any TPM elementary volume $dv$ can be thought of as delimited by the porous solid; therefore, when subjected to deviatoric strains, the volume change of an elementary volume $dv$ can directly be expressed as a function of the seepage velocity.

Furthermore, the mixture as a whole can be seen as the superposition of its constituents. Therefore both the balance of mass (Eq. \ref{eq:conservation of mass for constituent alpha}) and the balance of linear momentum (Eq. \ref{eq:linear momentum for constituent alpha}) are first expressed for each constituent independently. The rate at which momentum is transmitted by constituent $\alpha$ to the other constituent is accounted for in  $\mathbf{p}^{\alpha}$. 

\begin{equation}
\frac{\partial\rho^{\alpha}}{\partial t}+div\left(\rho^{\alpha}\mathbf{v}^{\alpha}\right)=0  \qquad \alpha = \left\{ S,F \right\}
\label{eq:conservation of mass for constituent alpha}
\end{equation}

\begin{equation}
div\left(\mathbf{T}^{\alpha}\right)+\mathbf{p}^{\alpha}=\mathbf{0}   \qquad \alpha = \left\{ S,F \right\}
\label{eq:linear momentum for constituent alpha}
\end{equation}

It is worth noting here that the apparent stress $\mathbf{T}^{\alpha}$ has the unit of a force per unit \emph{mixture area}, not unit area of constituent $\alpha$.

The governing equations of both constituents are then superimposed in order to obtain the governing equation of the mixture as a whole. Using (Eq. \ref{eq:connect density to true density}) and the fact that the true density is constant (assumption of intrinsic incompressibility), the mass balance equation for the mixture can be written as:

\begin{equation}
\frac{\partial\left(n^{S} + n^{F}\right)}{\partial t} + div\left(n^{S}\mathbf{v}^{S} + n^{F}\mathbf{v}^{F}\right)=0
\end{equation}

Using the saturation condition (Eq. \ref{eq:saturation condition}) and the seepage velocity (Eq. \ref{eq:seepage velocity}) allows us to rewrite the mass balance equation of the mixture:

\begin{equation}
div\left(\mathbf{v}^{S}+\mathbf{w}\right)=0
\label{eq:conservation of mass}
\end{equation}

It is interesting to realise that as a result of the ``smearing'' introduced by the concept of volume fractions, the seepage velocity $\mathbf{w}$ represents a macro-level average of the fluid velocity with respect to the solid phase and not an actual fluid velocity at pore-level (\cite{Cowin1999}).

In a similar fashion, the linear momentum balance of the mixture is obtained by summing (Eq. \ref{eq:linear momentum for constituent alpha}) for both constituents:

\begin{equation}
div\left(\mathbf{T}^{F} + \mathbf{T}^{S} \right) + \mathbf{p}^{F} + \mathbf{p}^{S}= \mathbf{0}
\label{eq:momentum - partial1}
\end{equation}

The introduction of the volume fractions has the consequence of introducing one extra variable for each phase of the mixture, which implies that the ``closure problem'' is not fulfilled (\cite{Boer1998}) (i.e. a greater number of unknowns than equations available). In order to overcome this problem, the saturation condition (Eq. \ref{eq:saturation condition}) is imposed as a kinematic constraint onto the entropy inequality, introducing a Lagrange multiplier $p$ (see for example \cite{Markert2005} or \cite{deBoer2005} and references therein for complete derivation). As a direct consequence, the concept of effective stress is introduced: 

\begin{equation}
\mathbf{T}^{\alpha} = -n^{\alpha}p\mathbf{I} + \mathbf{T}_{E}^{\alpha}
\label{eq:stress}
\end{equation}

The quantities $\mathbf{T}_{E}^{\alpha}$ are the constituent effective stresses, which must be determined constitutively. Summing over constituents, the mixture stress $\mathbf{T}$ is expressed in terms of the mixture effective stress $\mathbf{T}_{E}$:

\begin{subequations} 
        \begin{align} 
           &\mathbf{T}=\mathbf{T}^{S}+\mathbf{T}^{F}=-p\mathbf{I}+\mathbf{T}_{E} \label{subeqn:total stress} \\ 
           &\mathbf{T}_{E}=\mathbf{T}_{E}^{\text{S}}+\mathbf{T}_{E}^{\text{F}} \label{subeqn:effective stress} 
        \end{align} 
\label{eq:total and effective stress}
\end{subequations}

It is common (see for example discussion in \cite{Ehlers2008}) to consider the fluid as ideal (i.e. neglecting viscosity for slow draining materials such as soft tissues or hydrogels) and neglect the dissipative stress $\mathbf{T}^{F}_{E}$ in comparison to the drag generated by fluid-solid interactions $\mathbf{p}^{F}$, which is accounted for with Darcy's law. Hence, introducing $\bsig$ as a new notation for the solid effective stress and identifying $p$ as the fluid pressure, (Eq. \ref{subeqn:total stress}) simplifies to:

\begin{equation}
\mathbf{T} = -p\mathbf{I}+\bsig
\label{eq:partial stress}
\end{equation}

Truesdell's third metaphysical principle implies that the sum over all constituents of the partial balance relations has to take the same form as the balance of the single-phase material (i.e. in the current application $div \,\mathbf{T} =\mathbf{0}$). As a consequence, the constituent momentum exchanges cancel out: $\mathbf{p}^{F} + \mathbf{p}^{S}= \mathbf{0}$ and the linear momentum equation of the mixture is finally obtained by substituting (Eq. \ref{eq:partial stress}) into (Eq. \ref{eq:momentum - partial1}):

\begin{equation}
div \left( \bsig -p\mathbf{I} \right) = \mathbf{0}
\label{eq:linear momentum}
\end{equation}

The weak formulation of the problem is derived in a standard fashion (see Eq. \ref{eq:weak form} for expression in the current configuration). $\mathbf{f}$ and $g$ are the weighting functions of the linear momentum and mass balance respectively. In (Eq. \ref{subeqn:weak form mass}) $\mathbf{n}$ represents the outward normal vector and in (Eq. \ref{subeqn:weak form momentum}) $\mathbf{t}$, is the surface traction vector.

\begin{subequations} 
        \begin{align} 
           \int_{v} \left( g \, div \left( \mathbf{v} \right) - \mathbf{w}  \, \mathbf{\nabla}g \right) \, dv = \int_{a} g \, \mathbf{w} \, \mathbf{n} \, da \label{subeqn:weak form mass} \\ 
           \int_{v} \left( \nabla \mathbf{f} \right) : \left( \bsig - p\mathbf{I} \right) \, dv = \int_{a}  \mathbf{f} . \mathbf{t} \, da \label{subeqn:weak form momentum} 
        \end{align} 
\label{eq:weak form}
\end{subequations}

Finally, the standard Galerkin formulation is obtained by discretising the weak form in space using Taylor-Hood elements, where a linear approximations for the pressure field (four-node tetrahedron) and a quadratic one for the displacement field (ten-node tetrahedron) are used. The weighting functions are discretised using the same shape functions as those used for displacement and pressure. The discretised quantities can be expressed in matrix form as:
\begin{align}
           p &= \mathbf{N}^{p}\mathbf{p}^{e}  \qquad  \mathbf{u} = \mathbf{N}^{u}\mathbf{u}^{e}&  \nonumber \\ 
           g &= \mathbf{N}^{p} \mathbf{g}^{e} \qquad  \mathbf{f} = \mathbf{N}^{u}\mathbf{f^{e}}&
\label{eqref:discretization}
\end{align}

The variational form of (Eq. \ref{eq:weak form}) is:
\begin{subequations}
\allowdisplaybreaks
	\begin{align}
		R_{p} = \int_{v} \left\{ \Np \mathbf{g}^{e} \right\}^{\textrm{T}} \nabla  \mathbf{N}^{u} \dot{\mathbf{u}}^{e}  \, dv + \int_{v} k \left\{ \gNp \mathbf{g}^{e} \right\}^{\textrm{T}} \left\{ \gNp \mathbf{p}^{e} \right\} \, dv - \nonumber \\  \int_{a}  \left\{ \Np \mathbf{g}^{e} \right\}^{\textrm{T}} \bar{\mathbf{q}} \, da = 0&  \\
		R_{u} = \int_{v} \left\{ \gNu \mathbf{f}^{e} \right\}^{\textrm{T}}  \left( \bsig^{e}  +  \Np \mathbf{p}^{e}\right )\, dv - \int_{a}  \left\{ \Nu \mathbf{f}^{e} \right\}^{\textrm{T}} \bar{\mathbf{t}} \, da = 0 \\
		\quad R = R_{u} + \Delta t R_{p} = 0 \label{subeqn:compact form}
	\end{align}
	\label{eq:variational form}
\end{subequations}

A backward finite difference scheme is used for time integration (i.e. $\mathbf{\dot{u}_{t+\Delta t}}=(\mathbf{u}_{t+\Delta t}-\mathbf{u}_{t})/\Delta t$, where $\Delta t$ is the time increment). Introducing Darcy's law  with a permeability $k$ assumed constant for simplicity, $\mathbf{w} = -k \,\nabla p$, the linearized formulation in matrix form is derived:

\begin{equation}
	\left[
		\begin{array}{cc}
			\mathbf{K}_{uu} & \mathbf{K}_{up}\\
			\mathbf{K}_{up}^{\textrm{T}} & \Delta t \, \mathbf{K}_{pp}
		\end{array}
	\right]
	\left\{
		\begin{array}{c}
			\delta \mathbf{u} \\
			\delta \mathbf{p} 
		\end{array}
	\right\} =
	\left\{ 
		\begin{array}{c}
			\mathbf{F}_{u}^{ext}\\
			\Delta t \, \mathbf{F}_{p}^{ext}
		\end{array}
	\right\} -
	\left\{ 
		\begin{array}{c}
			\mathbf{F}_{u}^{int}\\
			\Delta t \, \mathbf{F}_{p}^{int}
		\end{array}
	\right\} 
\label{eq:matrix form}
\end{equation}

Each quantity $\left( \bullet \right)$ is linearised using the notations: $\left( \bullet \right) = \left( \hat{\bullet} \right) + \delta \left( \bullet \right)$. Terms in (Eq. \ref{eq:matrix form}) are detailed here:

\begin{subequations} 
\allowdisplaybreaks
        \begin{align} 
           &\mathbf{K}_{uu} =  \int_{v}  \gNut \mathbf{D} \gNu + \gNut \bsig \gNu \,dv\\ 
           &\mathbf{K}_{up} =  \int_{v} \Npt  \gNu\,dv\\ 
           &\mathbf{K}_{pp} =  \int_{v} \gNpt k \gNp\,dv\\ 
           &\mathbf{F}_{u}^{int} =  \int_{v} \gNut \hat{\bsig} + \gNut \Np \hatpe  \,dv\\
           &\mathbf{F}_{p}^{int} =  \int_{v} \gNpt k \gNp \hatpe + \Np \gNu \hat{\mathbf{v}}^{e} \,dv \\
           &\mathbf{F}_{u}^{ext} =  \int_{v}  \Nut \bar{\mathbf{t}} \,da \\
           &\mathbf{F}_{p}^{ext} =  \int_{v}  \Npt k \gNp \hatpe \,da
        \end{align} 
\label{eq:K and F terms}
\end{subequations}

where $\mathbf{D}$ is the matrix of material constants. In this work, the solid constituent is modelled with a Neo-Hooke model (see \cite{Bonet2008}):

\begin{equation}
\Phi_{NeoHooke} = \frac{\mu}{2} \left( I_{C} - 3 \right) - \mu \, ln\,J + \frac{\lambda}{2} \left( ln \, J \right)^{2} 
\label{eq:neohooke}
\end{equation}

The formulation was implemented into an in-house finite element code. The observed aforementioned instabilities (presented in the \emph{Numerical example} section) manifest in the form of spurious and non-physical oscillations in the pressure field in regions close to the free-flow boundaries. In the following section, we propose a stabilisation method to eliminate this issue.

\section{The GLS stabilization}

The spurious oscillations are stabilised using a Galerkin Least-Square (GLS) formulation, following \cite{Truty2001}. This formulation was originally derived to solve geotechnical problems of fully and partially saturated soils. Although the derivation is extended here for finite deformations, the outcome is simpler owing to the fact that the permeability is assumed constant. A weighted least-square term $R^{GLS}$, originating from the strong form of the fluid flow continuity equation (Eq. \ref{eq:conservation of mass}) is derived and added onto the weak form (Eq. \ref{subeqn:compact form}):

\begin{equation}
R = R_{u} - \Delta t R_{p} + R^{GLS} = 0
\end{equation}

Starting off by defining the least-square term $R^{GLS}$, where $\tau^{*}$ is a stabilisation factor that will subsequently be defined:

\begin{equation}
R^{GLS} =  \int_{v}  
\left[ div \left(\dot{\mathbf{f}} + k \nabla g \right)\right]^{\mathrm{T}}
\tau^{*}
\left[ div \left(\mathbf{v}^{S} + k \nabla p \right) \right] \, dv
\label{eq:GLS weak form}
\end{equation}

(Eq. \ref{eq:GLS weak form}) is rewritten as follows, taking into account the discretisation and the fact that the permeability is assumed constant and the pressure is linear:

\begin{equation}
R^{GLS} =  \int_{v}  
\left\{ \gNu \mathbf{\dot{f}}^{e} \right\}^{\textrm{T}}
\tau^{*}
\left\{ \gNu \mathbf{\dot{u}}^{e} \right\} \, dv
\label{eq:GLS discretised}
\end{equation}

Introducing the time integration scheme and defining $\tau^{GLS} = \tau^{*} / \left( \Delta t \right) ^{2}$:

\begin{equation}
R^{GLS} =  \int_{v}  
\left\{ \mathbf{N}^{u} \mathbf{f}^{e} \right\}^{\textrm{T}}
\tau^{GLS}
\left\{ \mathbf{N}^{u} \mathbf{u}^{e} \right\} \, dv
\label{eq:GLS discretised final}
\end{equation}

(Eq. \ref{eq:GLS discretised final}) is linearised using the notations $\mathbf{u}^{e} = \mathbf{\hat{u}} + \delta \mathbf{u}$ and writes in matrix form:

\begin{equation}
	\left[
		\begin{array}{cc}
			\mathbf{K}_{uu} + \mathbf{K}^{GLS}_{uu} & \mathbf{K}_{up}\\
			\mathbf{K}_{up}^{\textrm{T}} & \Delta t \, \mathbf{K}_{pp}
		\end{array}
	\right]
	\left\{
		\begin{array}{c}
			\delta \mathbf{u} \\
			\delta \mathbf{p} 
		\end{array}
	\right\} =
	\left\{ 
		\begin{array}{c}
			\mathbf{F}_{u}^{ext}\\
			\Delta t \, \mathbf{F}_{p}^{ext}
		\end{array}
	\right\} -
	\left\{ 
		\begin{array}{c}
			\mathbf{F}_{u}^{int} + \mathbf{F}^{GLS-int}_{u}\\
			\Delta t \, \mathbf{F}_{p}^{int}
		\end{array}
	\right\} 
\label{eq:GLS matrix form}
\end{equation}

where:

\begin{subequations} 
        \begin{align} 
           \mathbf{K}^{GLS}_{uu} = \int_{v} \left\{ \gNu \right\}^{\textrm{T}} \tau^{GLS} \left\{ \gNu \right\} \, dv\label{subeqn:K_GLS} \\ 
           \mathbf{F}^{GLS-int}_{u} = \int_{v} \left\{ \gNu \right\}^{\textrm{T}} \tau^{GLS} \left\{ \gNu \mathbf{\hat{u}}^{e} \right\} \, dv \label{subeqn:FintU_GLS} 
        \end{align} 
\label{eq:GLS terms}
\end{subequations}

Finally, the stabilization factor is defined based on \cite{Truty2001}:

\begin{equation}
\tau^{GLS} = \frac{h^2}{4k\Delta t}
\label{eq:tau GLS}
\end{equation}

In \cite{Truty2001}, the parameter $h$ is a length characterizing the element's size in the direction of the fluid flow. In the current application, it was simply taken as the radius of the element's circumsphere. As illustrated in the following section, this definition of $\tau^{GLS}$ stems from the fact that the solution needs greater stabilisation as the mesh gets coarser and as either or both the permeability and the size of the time-step decrease. In this work, due to the soft nature of the solid phase, it was not necessary to include a measure of material stiffness in the definition of $\tau^{GLS}$. Although this may be necessary for stiffer materials (e.g. \cite{Truty2001} and \cite{Aguilar2008}).

It can be noted that the numerical integration is performed using a 4-point gaussian quadrature. The accuracy was verified with tests using up to 27 points. Reduced integration was also considered for the GLS terms, with no beneficial outcome.

\section{Numerical example}

The performance of the GLS stabilisation technique is assessed on a biphasic cylinder subjected to unconfined compression. With a 18mm radius, a thickness of 8mm and a solid phase defined with $\lambda$ = 0.2 MPa and $\mu$ = 0.5 MPa, the cylinder can be thought of as an idealised human nucleus pulposus of the intervertebral disc when the permeability is set to $k$ = 1$\times$10$^{-3} $mm$^{4}$N$^{-1}$s$^{-1}$.

\begin{figure}[ht]
\centering
\includegraphics[width=0.65\textwidth]{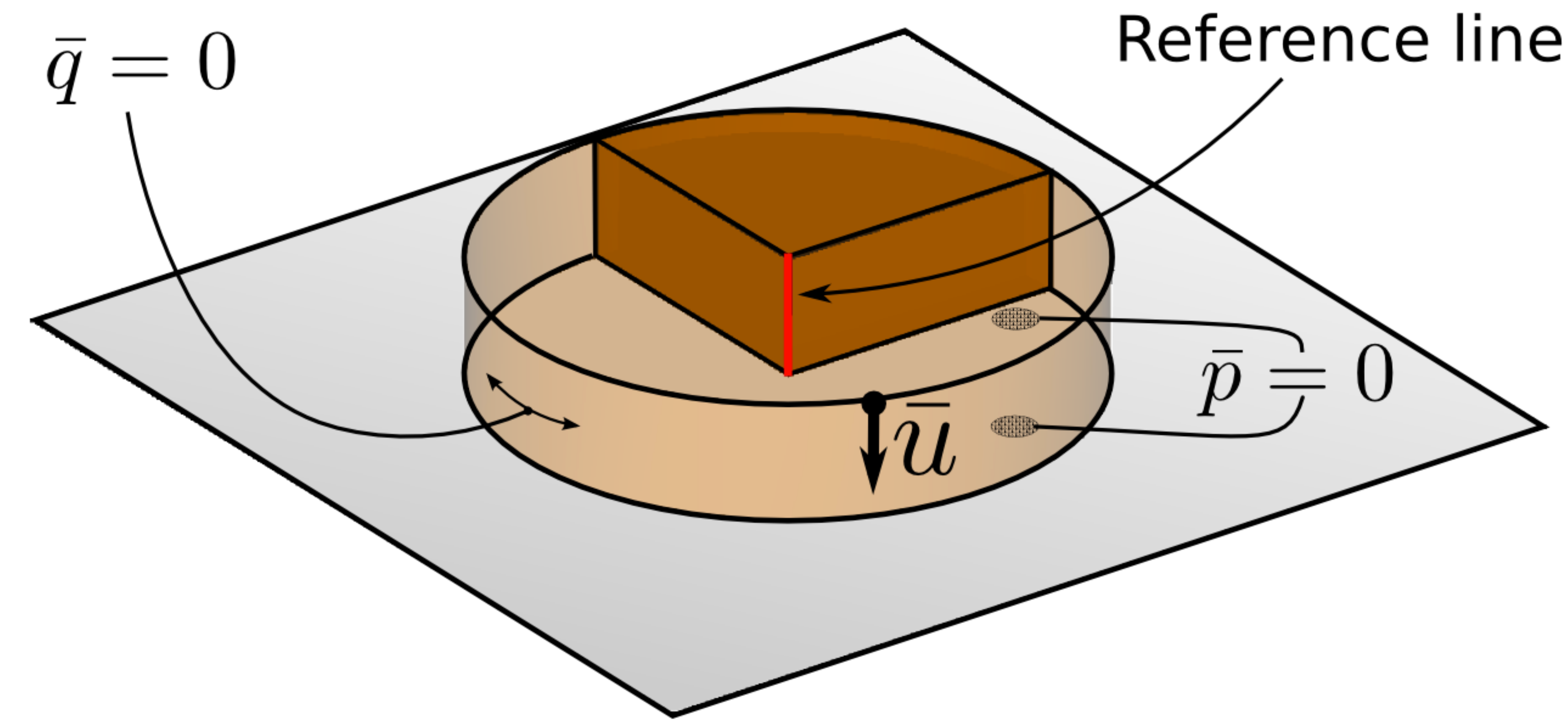}
\caption{Loading and boundary conditions}
\label{fig:Loading and boundary conditions}
\end{figure}

The fluid flux $\bar{\mathbf{q}}$ at the boundary is prescribed to zero on the vertical faces offering a lateral seal to the cylinder. The pressure  $\bar{p}$ is set to zero on the top and bottom surfaces (see Fig. \ref{fig:Loading and boundary conditions}). These boundary conditions define the top and bottom surfaces as the only free-flow boundaries. This results in a near to uni-axial fluid flow through the depth of the cylinder.

\begin{figure}[ht]
\begin{minipage}[b]{0.3\linewidth}
\centering
	\begin{tabular}{c|ccc}
          \hline
	  & \#		&  \#	 	& \# elements \\
	  & Nodes	&  Elements	& on ref line \\
          \hline
   mesh 1 & 8083	& 5261	 	& 5	\\
   mesh 2 & 16507	& 11087	 	& 7	\\
   mesh 3 & 38424	& 26622	 	& 9	\\
   mesh 4 & 78033	& 55190	 	& 12	\\
          \hline
	\end{tabular}
\par\vspace{0pt}
\end{minipage}
\begin{minipage}[t]{1.\linewidth}
\centering
    \includegraphics[scale = 0.2]{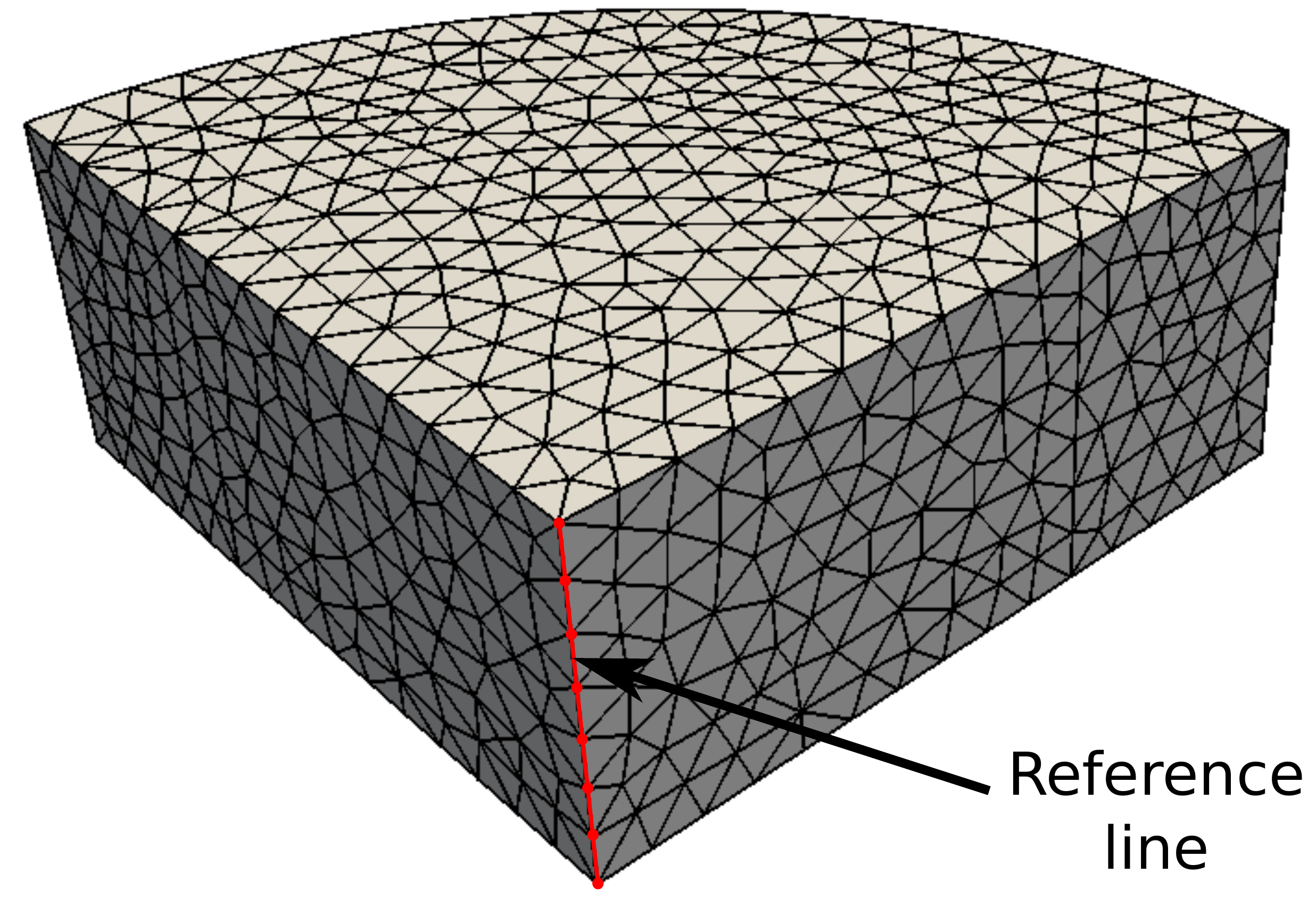}
\par\vspace{0pt}
\end{minipage}
\caption{Meshes characteristics (left) and mesh 2 (right)}
\label{fig:mesh description}
\end{figure}

In what is herein referred to as the ``reference loading'', the top surface is displaced downwards at a rate of 2.5$\mu m.s^{-1}$ with time increments $\Delta t$ = 6.4s, until the height of the cylinder reduces by 1\%. Such loading rates, together with permeabilities lower than 1 mm$^{4}$ N$^{-1}$ s$^{-1}$, guaranties that the steady state is not reached instantly. In order to reduce the size of the problem, symmetry boundary conditions are applied onto a quarter cylinder. Analyses are undertaken on different meshes, the main characteristics of which are shown in Fig. \ref{fig:mesh description}.

All results are plotted against nodal values gathered on the ``reference line'' defined in Fig. \ref{fig:Loading and boundary conditions} and \ref{fig:mesh description}.

\begin{figure}[H]
  \subfloat[Nodal strains]{\label{fig:mesh2 - k0 - strain}\includegraphics[width=0.33\textwidth]{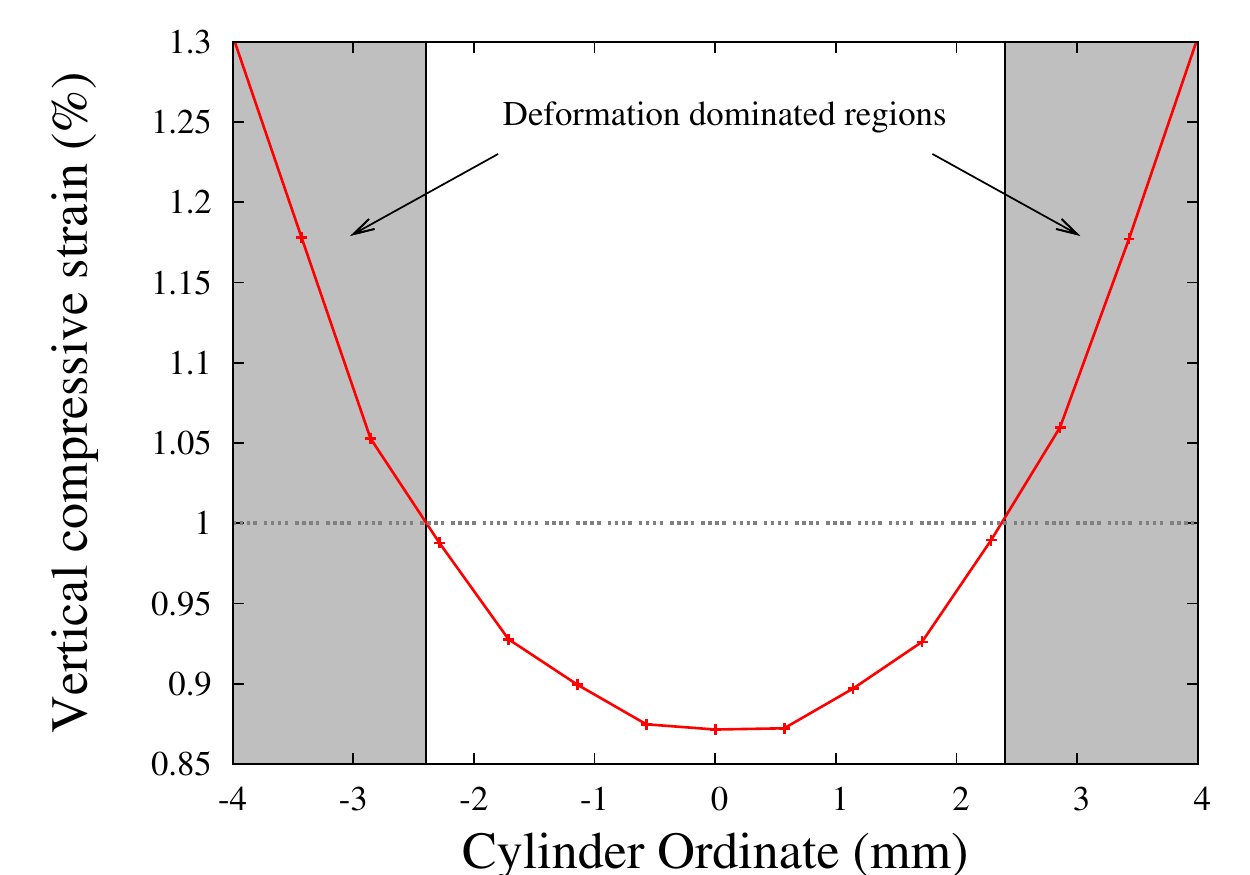}}
  \subfloat[Nodal pressures]{\label{fig:mesh2 - k0 - press}\includegraphics[width=0.33\textwidth]{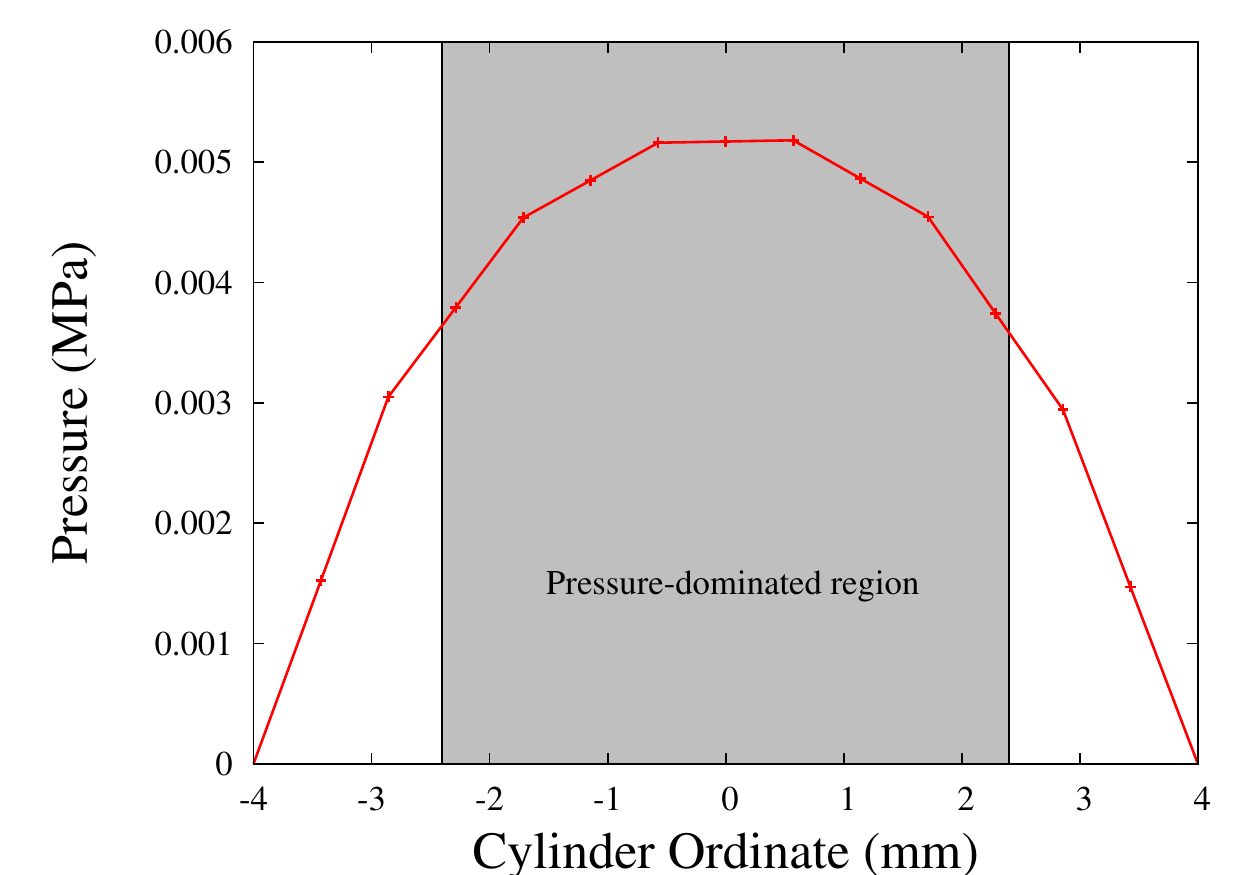}} 
  \subfloat[Influence of permeability]{\label{fig:pb mesh2 noGLS vary perm}\includegraphics[width=0.33\textwidth]{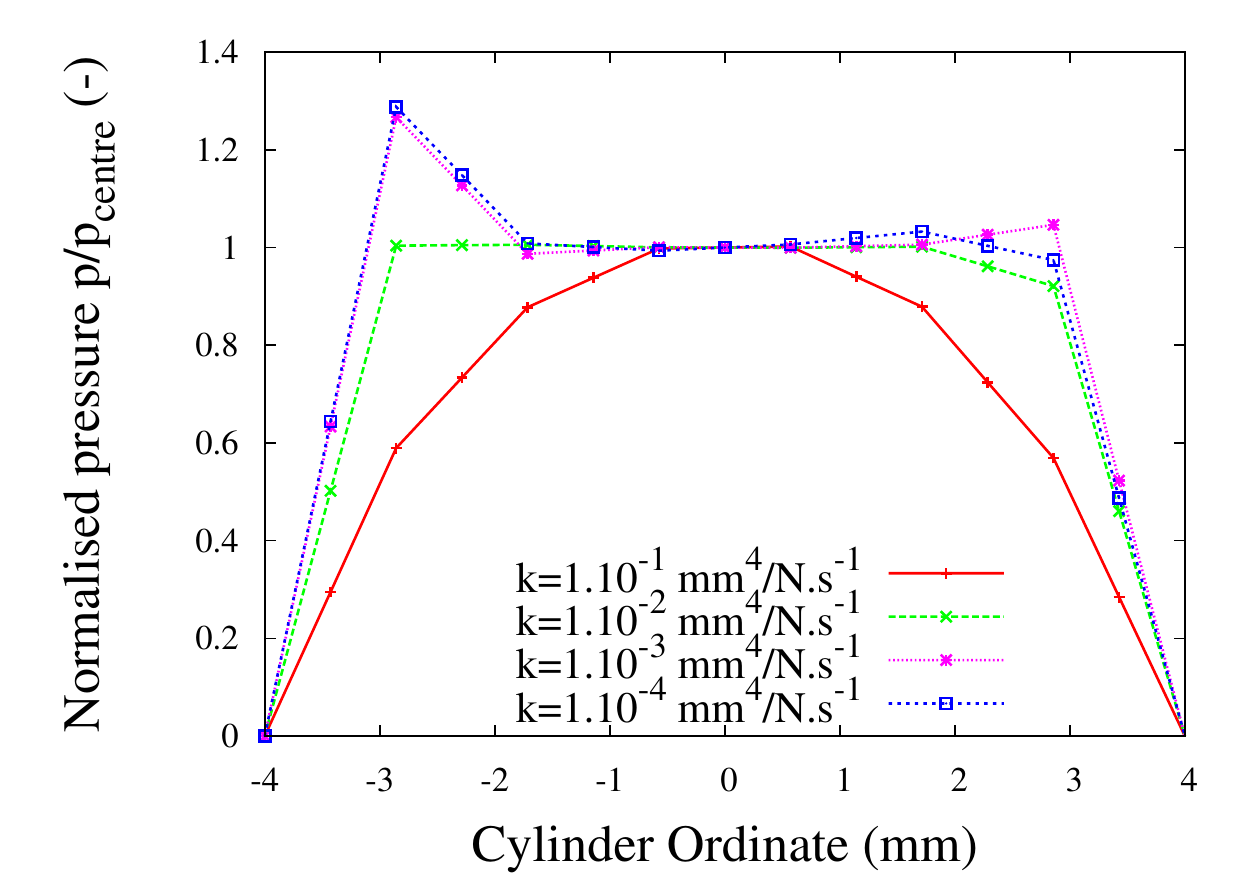}}
 \caption{Mesh 2: deformation mechanisms for $k$ = 0.1 mm$^{4}$N$^{-1}$s$^{-1}$ (left and centre) and influence of permeability (right)}
 \label{fig:expected deformations} 
\end{figure}

The low permeability hinders the fluid's ability to flow, defining two distinct regions associated with the load transfer mechanism. The first region, located near the free-flow boundaries, is dominated by solid deformation: the fluid does not have the ability to pressurize as it is squeezed out of the cylinder. This results in a lower fluid content in this region, explaining the peak strains observed near the top and bottom surfaces (see Fig. \ref{fig:mesh2 - k0 - strain}). The second region, situated at the centre of the cylinder, is predominantly subjected to fluid pressurization (see Fig. \ref{fig:mesh2 - k0 - press}) due to the fact that the low permeability is confining the fluid at the centre of the cylinder.

\begin{figure}[H]
\subfloat[$k$ = 1$\times$10$^{-2}$ mm$^{4}$N$^{-1}$s$^{-1}$]{\label{fig:mesh - k2 - noGLS}\includegraphics[width=0.33\textwidth]{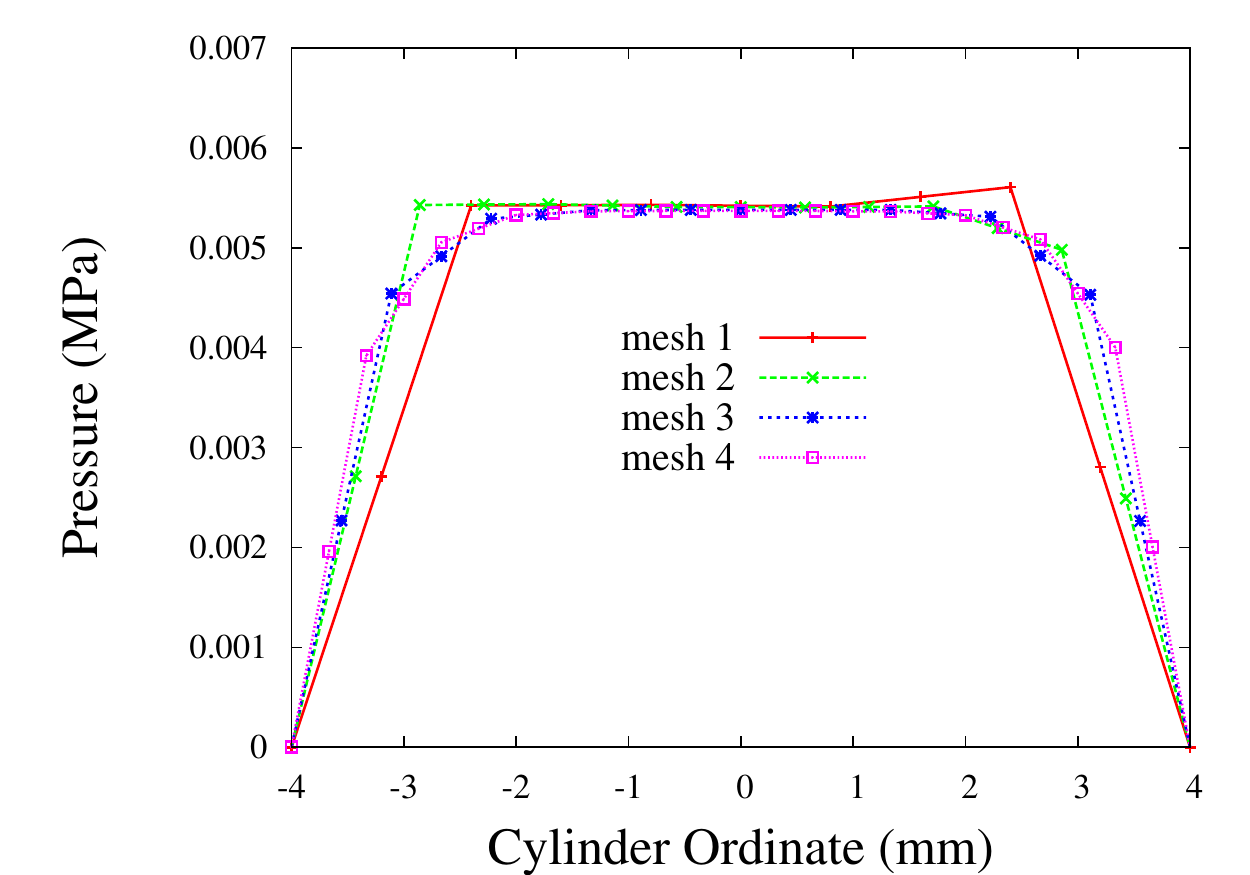}}
  \subfloat[$k$ = 5$\times$10$^{-3}$ mm$^{4}$N$^{-1}$s$^{-1}$]{\label{fig:mesh - 5k3 - noGLS}\includegraphics[width=0.33\textwidth]{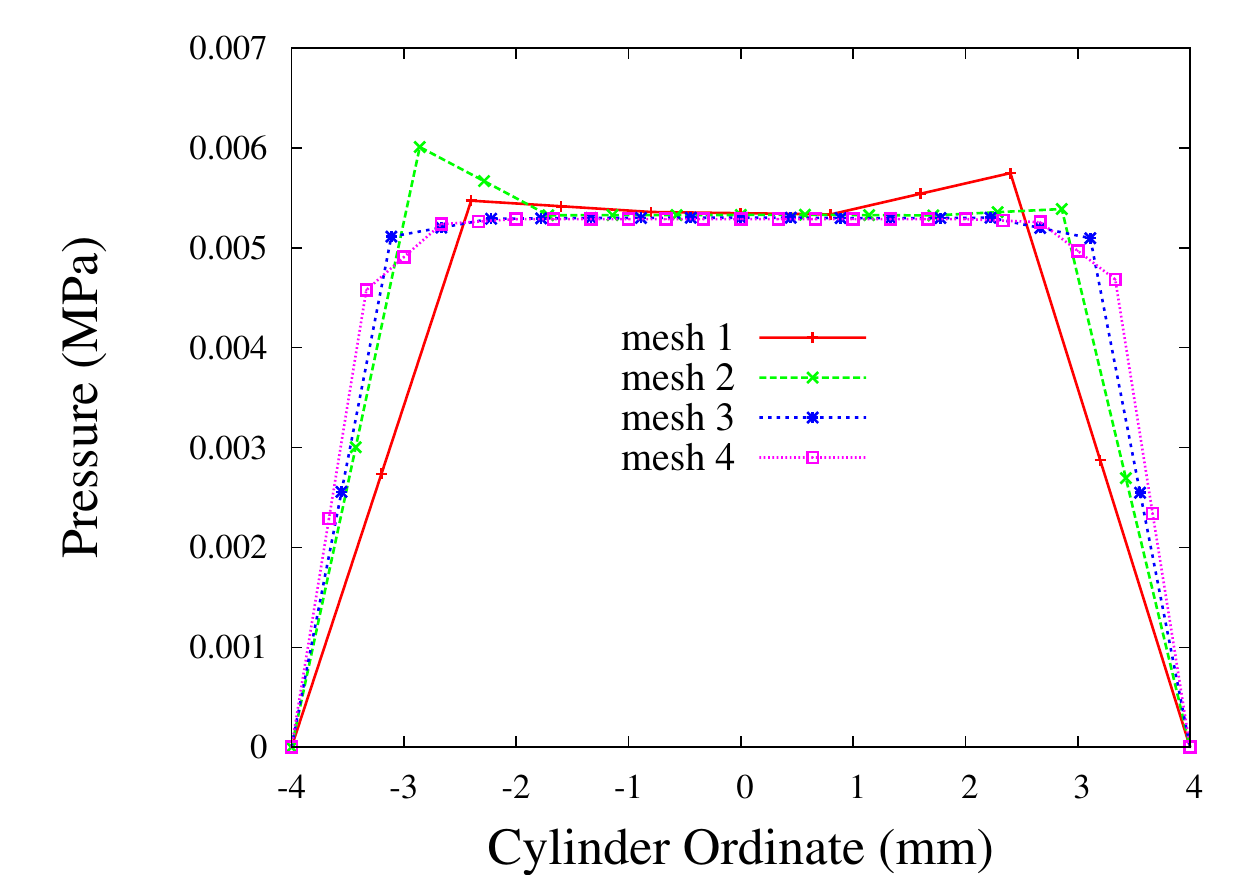}} 
  \subfloat[$k$ = 1$\times$10$^{-3}$ mm$^{4}$N$^{-1}$s$^{-1}$]{\label{fig:mesh - k3 - noGLS}\includegraphics[width=0.33\textwidth]{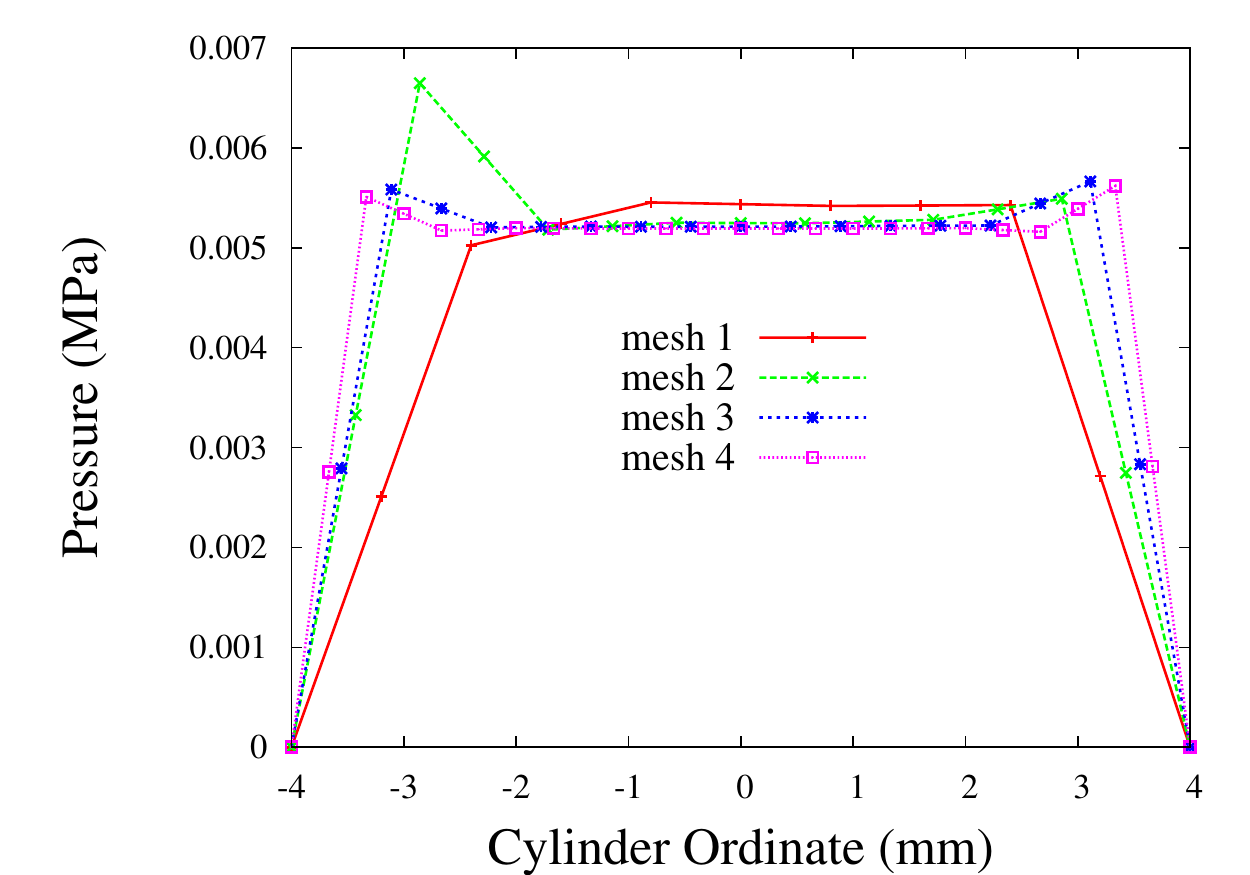}}
 \caption{Effects of mesh refinement on standard Galerkin method}
 \label{fig:Mesh refinement, no GLS} 
\end{figure}

As the permeability decreases, the boundary between solid- and pressure-dominated regions shifts towards the top and bottom surfaces (see Fig. \ref{fig:pb mesh2 noGLS vary perm}) and the level of pressurization rises. For a given mesh, when the permeability falls under a certain value (k $<$ 1$\times$10$^{-1}$ mm$^{4}$N$^{-1}$s$^{-1}$ for mesh 2 in this example), the pressure profile starts to exhibit spurious oscillations near the free-flow boundaries (over 10\% discrepancies for mesh 1 and mesh 2 when k = 5$\times$10$^{-2}$ mm$^{4}$N$^{-1}$s$^{-1}$, and over 8\% for all meshes when k = 1$\times$10$^{-2}$ mm$^{4}$N$^{-1}$s$^{-1}$). Furthermore, the quality of the solution can be affected through the entire mesh as the near boundary oscillation propagates toward the centre for coarse meshes (e.g. mesh 1 in Fig. \ref{fig:mesh - k3 - noGLS}). Finally, it was verified (in line with \cite{Vermeer1981}) that decreasing the time-step exaggerates the pressure oscillations, although not presented here.

Mesh refinement is the most natural and straightforward choice to overcome this issue, in particular when interested in accurately capturing the steep pressure gradients. As Fig. \ref{fig:mesh - k2 - noGLS} and \ref{fig:mesh - 5k3 - noGLS} illustrate, the non-physical pressure peaks can be removed by using meshes denser than that characterised by mesh 2. However, this is only a valid solution some cases. Fig. \ref{fig:mesh - k3 - noGLS} illustrates that the spurious oscillations cannot always be reduced by reasonably-sized denser meshes for low permeability.

\begin{figure}[H]
  \subfloat[$k$ = 5$\times$10$^{-3}$ mm$^{4}$N$^{-1}$s$^{-1}$]{\label{fig:5k3 GLS vary mesh}\includegraphics[width=0.5\textwidth]{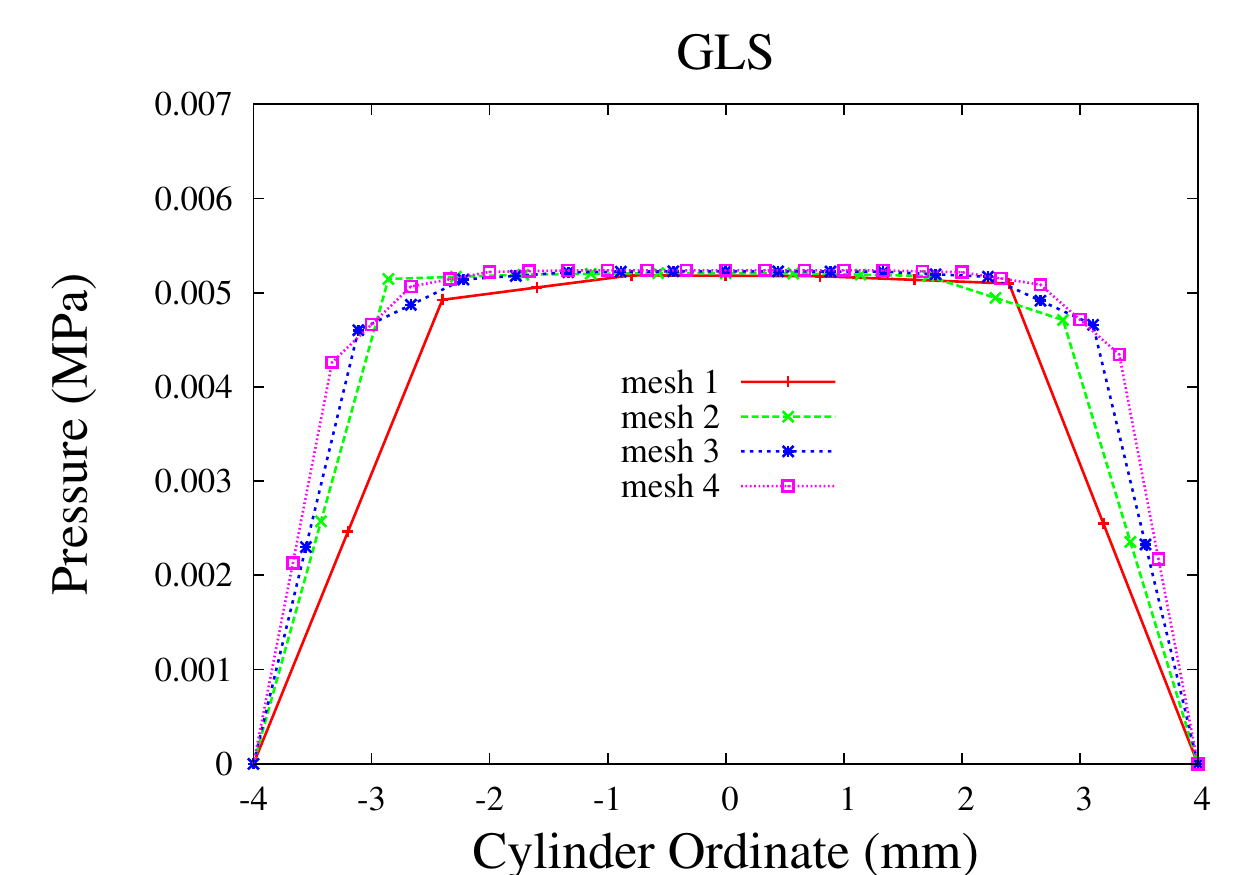}} 
  \subfloat[$k$ = 1$\times$10$^{-3}$ mm$^{4}$N$^{-1}$s$^{-1}$]{\label{fig:k3 GLS vary mesh}\includegraphics[width=0.5\textwidth]{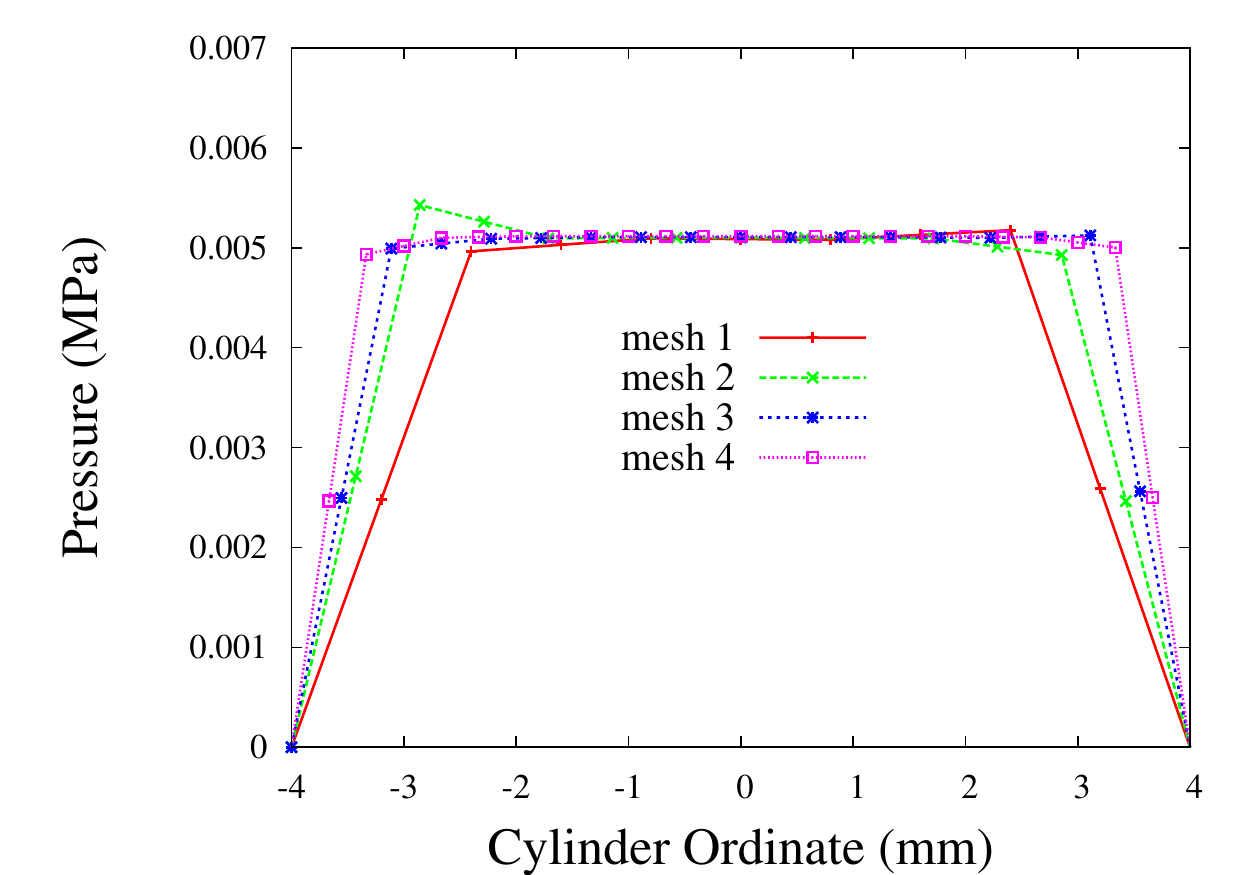}}
 \caption{Effects of mesh refinement on GLS stabilisation}
 \label{fig:Mesh refinement, GLS}
\end{figure}

The GLS stabilization offers substantial improvements to the solution. Fig. \ref{fig:Mesh refinement, GLS} shows the benefits for four different meshes when $k$ = 5$\times$10$^{-3}$ mm$^{4}$N$^{-1}$s$^{-1}$ and $k$ = 1$\times$10$^{-3}$ mm$^{4}$N$^{-1}$s$^{-1}$. The primary enhancement is that all spurious oscillations observed in Fig. \ref{fig:Mesh refinement, no GLS} have been stabilised, with the exception of mesh 2 where, when $k$ = 1$\times$10$^{-3}$ mm$^{4}$N$^{-1}$s$^{-1}$, the discrepancies decreased from 27\% to 6\%. Additionally, when the accurate resolution of the near-boundary pressure gradient, defining the transition between the deformation- and pressure-driven regions, is not sought, the GLS formulation allows for coarser meshes to be used, since it also prevents the oscillations from propagating towards the centre. For example, when $k$ = 1$\times$10$^{-3}$ mm$^{4}$N$^{-1}$s$^{-1}$, mesh 1 (8000 nodes) with GLS offers similar performances to mesh 3 (38000) without the stabilisation (compare Fig. \ref{fig:mesh - k3 - noGLS} and \ref{fig:k3 GLS vary mesh}). Finally, it is important to notice that the GLS stabilisation is only having a damping effect on the spurious oscillations, while leaving stable solutions unaffected (compare Fig. \ref{fig:mesh - 5k3 - noGLS} and \ref{fig:5k3 GLS vary mesh}).

A few observations can be made to support the choice of the stabilization factor $\tau^{GLS}$  in (Eq. \ref{eq:tau GLS}). First, Fig. \ref{fig:Mesh refinement, GLS} highlights the fact that the stabilization performs equally well irrespective of the change of mesh and permeability, giving confidence in the way the element's characteristic size and permeability are accounted for. The impact of the time-step on the stabilization was also investigated: it was verified that if $\tau^{GLS}$ is not inversely proportional to $\Delta t$, stabilisation is not possible. Also, in simulations not presented here, it was confirmed that the size of the time-step ($\Delta t$ = \{0.64s, 3.2s, 6.4s, 8s, 32s\}) does not affect the quality of the stabilised solution. Finally, it was also verified (again not shown here) that changing the loading rate (1.25$\mu m.s^{-1}$, 2.5$\mu m.s^{-1}$, 6.25$\mu m.s^{-1}$) did not affect the degree of peak pressure oscillations for the stabilised results.

\begin{figure}[H]
  \subfloat[Standard Galerkin]{\label{fig:q2 k3 increase strain noGLS}\includegraphics[width=0.5\textwidth]{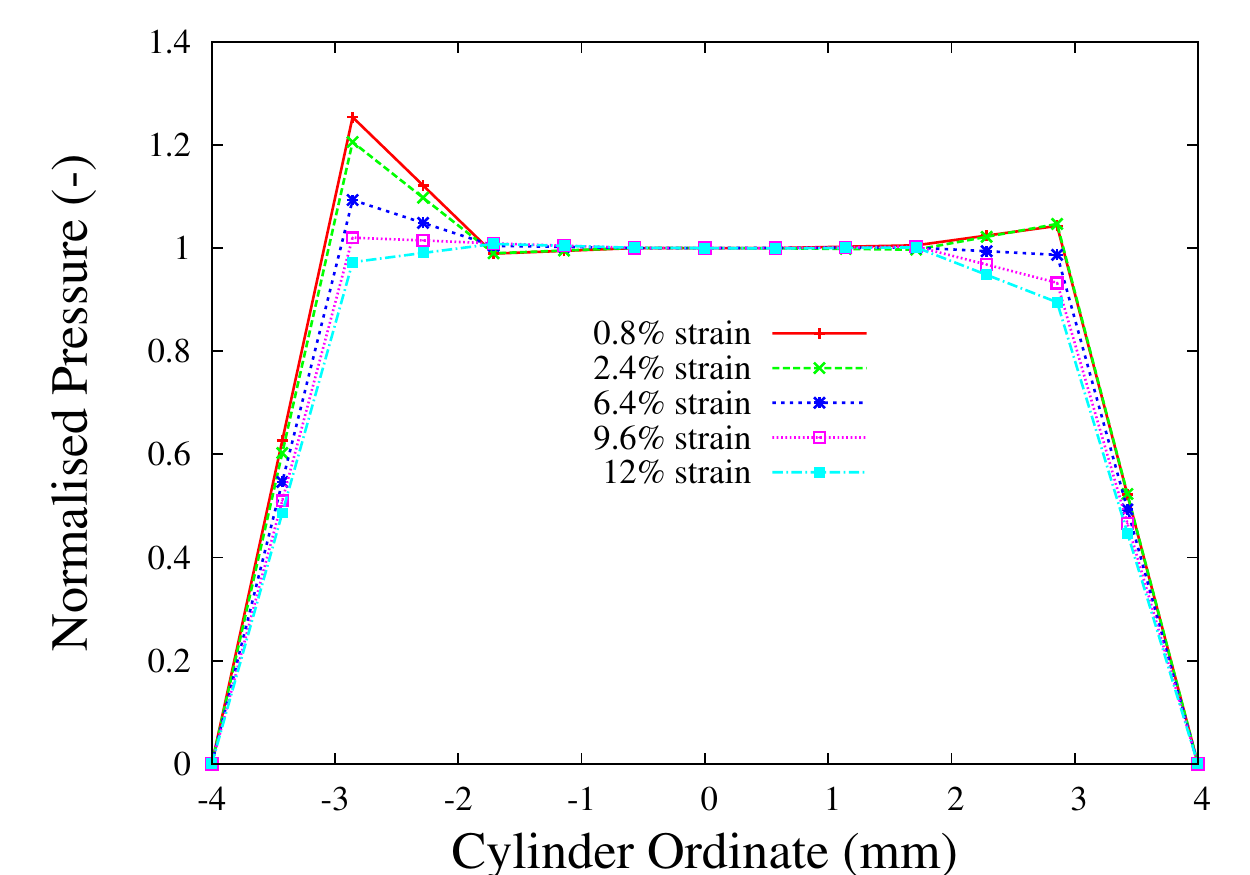}}
  \subfloat[GLS]{\label{fig:q2 k3 increase strain GLS}\includegraphics[width=0.5\textwidth]{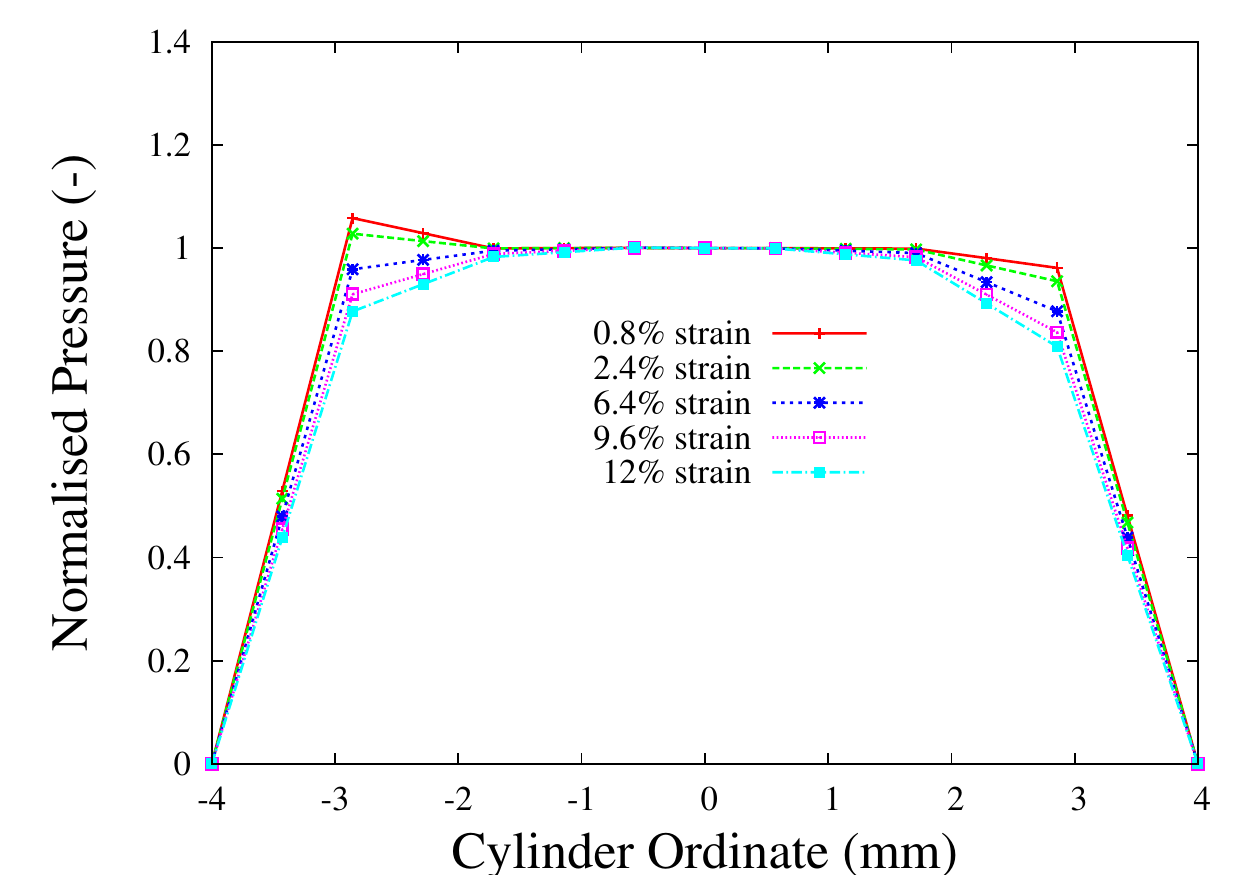}} 
 \caption{Effects of higher strain on stability for mesh 2}
 \label{fig:q2 k3 increase strain}
\end{figure}

Performance of the GLS stabilization was initially assessed for greater levels of deformation. As Fig. \ref{fig:q2 k3 increase strain} shows, the level of pressure discrepancies reduces as the compressive strain increases, which is in line with the findings in \cite{Truty2001} and \cite{Aguilar2008}, where oscillations are reported to occur at the ``early stage'' of the consolidation problem. Although the GLS stabilisation also performs well at higher strains (see Fig. \ref{fig:q2 k3 increase strain GLS}), this observation motivated the choice to present results at 1\% compression throughout this numerical example.

A sensitivity study was performed to charaterise the parameter $h$ in (Eq. \ref{eq:tau GLS}). Several combinations of the radius of the circumsphere, the shortest and longest edges of the tetrahedron were tested without noticeable and consistent improvement to the overall solution.

\section{Conclusion}

It was observed that an hyperelastic biphasic model, implemented in a finite element framework with Taylor-Hood tetrahedral elements, exhibits non-physical pressure oscillations for low permeabilities. A Galerkin least-square formulation was derived for finite deformations in order to stabilise these oscillations.

In the context of constant permeability and near to uni-axial fluid flow, the current formulation shows good results. It eliminates the spurious oscillations for most meshes (and damp the oscillations for others meshes) and also prevents these oscillations from propagating towards the centre of the medium as reported for very coarse meshes. The solution scheme proved to be robust when tested against various mesh densities, permeabilities, loading rates, compressive strains and time steps. It is also worth mentioning that the benefits of this formulation come at minimal computational cost, as no additional degrees of freedom are required.

For more complex fluid flow situations, information regarding the fluid flow directionality should be included in the derivation of the element's characteristic size. In \cite{Truty2001}, this has been considered for 2D and for 3D brick elements. It will prove more challenging for 3D tetrahedral elements. For extension to strain-dependent permeability, derivation of the stabilisation terms is straightforward but would result in a non-symmetrical system that would be computationally more expensive to solve.

\section{Acknowledgements}
This work was supported by the Glasgow Research Partnership in Engineering (GRPE).

\bibliography{JabRef-Database}
\end{document}